\documentclass[titlepage,draft,12pt]{article} 
\usepackage{amssymb,amsthm,amsmath}  
\usepackage[utf8]{inputenc}
\usepackage[a4paper]{geometry}

\date{}

\newcommand{\re}{\mathbb{R}}

\newcommand{\ep}{\varepsilon}
\renewcommand{\qed}{{\penalty 10000\mbox{$\quad\Box$}}}

\newcommand{\aep}{a_{\ep}}
\newcommand{\bep}{b_{\ep}}
\newcommand{\gep}{g_{\ep}}
\newcommand{\yep}{y_{\ep}}
\newcommand{\uep}{u_{\ep}}

\newtheorem{thm}{Theorem}[section]
\newtheorem{thmbibl}{Theorem}
\newtheorem{propbibl}[thmbibl]{Proposition}

\newtheorem{prop}[thm]{Proposition}

\newtheorem{lemma}[thm]{Lemma}

\title{Symmetry-breaking in a generalized Wirtinger inequality}

 \author{Marina Ghisi\vspace{1ex}\\ 
 {\normalsize Universit\`a degli Studi di Pisa} \\
 {\normalsize Dipartimento di Matematica}\\ 
 {\normalsize PISA (Italy)}\\
 {\normalsize e-mail: \texttt{marina.ghisi@unipi.it}}
 \and
 Massimo Gobbino\vspace{1ex}\\ 
 {\normalsize Universit\`a degli Studi di Pisa} \\
 {\normalsize Dipartimento di Ingegneria Civile e Industriale}\\ 
 {\normalsize PISA (Italy)}\\  
 {\normalsize e-mail: \texttt{massimo.gobbino@unipi.it}}
 \and
 Giulio Rovellini\vspace{1ex}\\ 
 {\normalsize Scuola Normale Superiore} \\
 {\normalsize Classe di Scienze}\\ 
 {\normalsize PISA (Italy)}\\  
 {\normalsize e-mail: \texttt{giulio.rovellini@sns.it}}}

\begin{document}
\maketitle
\begin{abstract}

The search of the optimal constant for a generalized Wirtinger inequality in an interval consists in minimizing the $p$-norm of the derivative among all functions whose $q$-norm is equal to~1 and whose $(r-1)$-power has zero average. Symmetry properties of minimizers have attracted great attention in mathematical literature in the last decades, leading to a precise characterization of symmetry and asymmetry regions.
	
In this paper we provide a proof of the symmetry result without computer assisted steps, and a proof of the asymmetry result which works as well for local minimizers. As a consequence, we have now a full elementary description of symmetry and asymmetry cases, both for global and for local minima. 

Proofs rely on appropriate nonlinear variable changes.

\vspace{6ex}

\noindent{\bf Mathematics Subject Classification 2010 (MSC2010):}
26D10, 49R05.

% 26D10   	Inequalities involving derivatives and differential and integral operators
% 49R05   	Variational methods for eigenvalues of operators

\vspace{6ex}


\vspace{6ex}

% \listoftodos

\end{abstract}

%%%%%%%%%%%%%%%%%%%%%
%                   %
%   Inizio lavoro   %
%                   %
%%%%%%%%%%%%%%%%%%%%%
 
\section{Introduction}

Let $(a,b)\subseteq\re$ be an interval. It is well-known that the inequality
\begin{equation}
\pi^{2}\int_{a}^{b}|u(x)|^{2}\,dx\leq(b-a)^{2}\int_{a}^{b}|u'(x)|^{2}\,dx
\nonumber
\end{equation}
holds true\begin{itemize}
  \item (\emph{Poincaré inequality}) for every $u\in C^{1}([a,b])$ such that $u(a)=u(b)=0$,
  
  \item (\emph{Wirtinger inequality}) for every $u\in C^{1}([a,b])$ such that
\begin{equation}
\int_{a}^{b}u(x)\,dx=0.
\nonumber
\end{equation}

\end{itemize}

More generally, given three real numbers $p>1$, $q>1$, $r>1$, after defining $p_{*}$ and $\theta$ in such a way that
\begin{equation}
\frac{1}{p}+\frac{1}{p_{*}}=1
\qquad\mbox{and}\qquad
\theta:=\frac{1}{p_{*}}+\frac{1}{q},
\nonumber
\end{equation}
one can show the existence of
\begin{itemize}

\item  (\emph{generalized Poincaré inequality}) a largest positive constant $\lambda_{P}(p,q)$ such that
\begin{equation}
\lambda_{P}(p,q)\left(\int_{a}^{b}|u(x)|^{q}\,dx\right)^{1/q}\leq(b-a)^{\theta}\left(\int_{a}^{b}|u'(x)|^{p}\,dx\right)^{1/p}
\label{ineq:P}
\end{equation}
for every $u\in C^{1}([a,b])$ such that $u(a)=u(b)=0$,

\item  (\emph{generalized Wirtinger inequality}) a largest positive constant $\lambda_{W}(p,q,r)$ such that
\begin{equation}
\lambda_{W}(p,q,r)\left(\int_{a}^{b}|u(x)|^{q}\,dx\right)^{1/q}\leq(b-a)^{\theta}\left(\int_{a}^{b}|u'(x)|^{p}\,dx\right)^{1/p}
\label{ineq:W}
\end{equation}
for every $u\in C^{1}([a,b])$ such that
\begin{equation}
\int_{a}^{b}|u(x)|^{r-2}u(x)\,dx=0.
\label{hp:W-r}
\end{equation}

\end{itemize}

Several problems in the calculus of variations can be reduced to $\lambda_{W}(p,q,r)$, some of which are listed in the introduction to~\cite{nazarov-2}. We also refer to~\cite{dacorogna-2} for further motivations of this search for optimal constants in Sobolev type inequalities.

The optimal constants in the two inequalities above can be characterized as the minimum of the quotient
\begin{equation}
(b-a)^{\theta}\frac{\|u'\|_{L^{p}((a,b))}}{\|u\|_{L^{q}((a,b))}}
\label{quotient}
\end{equation}
among all intervals $[a,b]\subseteq\re$ and all functions $u\in C^{1}([a,b])$ that do not vanish identically in $[a,b]$ and satisfy
\begin{itemize}
\item  the boundary conditions $u(a)=u(b)=0$ in the case of Poincaré inequality,

\item  the integral condition (\ref{hp:W-r}) in the case of Wirtinger inequality.
\end{itemize}

The existence of $\lambda_{P}(p,q)$ and $\lambda_{W}(p,q,r)$ can be proved by applying in a standard way the direct method in the calculus of variations (see~\cite{dacorogna-1,dacorogna-2}). As a by-product, inequalities (\ref{ineq:P}) and (\ref{ineq:W}) hold true with the same constants even if we replace $u\in C^{1}([a,b])$ by $u\in W^{1,p}((a,b))$, of course subject to the same boundary or integral conditions. 

Due to scale invariance, there is no loss of generality in working in a fixed interval, for example $[-1,1]$.  Focussing on this interval, in the special case $p=q=r=2$ one can compute explicitly the optimal constants and characterize the equality cases as follows:
\begin{itemize}
\item  $\lambda_{W}(2,2,2)=\lambda_{P}(2,2)=\pi$, 

\item  minimizers (both local and global) to $\lambda_{P}(2,2)$ are even functions (actually all nonzero multiples of $\cos(\frac{\pi}{2}x)$), and they are eigenvectors relative to the first eigenvalue of the 1-d Laplacian with Dirichlet boundary conditions,

\item  minimizers (both local and global) to $\lambda_{W}(2,2,2)$ are odd functions (actually all nonzero multiples of $\sin(\frac{\pi}{2}x)$),  and they are eigenvectors relative to the second eigenvalue of the 1-d Laplacian with Neumann boundary conditions,

\item  minimizers to $\lambda_{W}(2,2,2)$ can be obtained from minimizers to $\lambda_{P}(2,2)$ with a ``cut-and-paste'' procedure, in the sense that if $u_{P}(x)$ is any equality case for Poincaré inequality, then
\begin{equation}
u_{W}(x):=\left\{\begin{array}{ll}
-u_{P}(x+1) & \mbox{if }x\in[-1,0]  \\
\noalign{\vspace{1ex}}
u_{P}(x-1) & \mbox{if }x\in[0,1] 
\end{array}\right.
\label{cut-and-paste}
\end{equation}
realizes the equality in Wirtinger inequality.
\end{itemize}

It is reasonable to ask whether these symmetry properties remain true for general values of the parameters $p$, $q$, $r$. This problem has generated a lot of literature, leading to the following answer.

\begin{thmbibl}[Symmetry/asymmetry in generalized Poincaré/Wirtinger inequalities]\label{thmbibl:main}

Let $p>1$, $q>1$, $r>1$ be three real numbers.

Then the following statements hold true.
\begin{enumerate}
\renewcommand{\labelenumi}{(\arabic{enumi})}
\item \emph{(Poincaré inequality -- Symmetry of minimizers)} Minimizers (both local and global) to $\lambda_{P}(p,q)$ are  always even functions, and $\lambda_{W}(p,q,r)\leq\lambda_{P}(p,q)$ for every admissible choice of the parameters. 

\item  \emph{(Wirtinger inequality -- Symmetry of minimizers)} If $q\leq(2r-1)p$ it turns out that $\lambda_{W}(p,q,r)=\lambda_{P}(p,q)$. Moreover, minimizers (both local and global) to $\lambda_{W}(p,q,r)$ are the odd functions obtained from minimizer to $\lambda_{P}(p,q)$ through the cut-and-paste procedure (\ref{cut-and-paste}).

\item  \emph{(Wirtinger inequality -- Asymmetry of minimizers)} If $q>(2r-1)p$ it turns out that $\lambda_{W}(p,q,r)<\lambda_{P}(p,q)$, and no odd function can be a minimizer for $\lambda_{W}(p,q,r)$ (not even a local minimizer).
\end{enumerate}

\end{thmbibl}

The proof of most parts of Theorem~\ref{thmbibl:main} was achieved in a series of papers of the last 25~years. The following table sums up the main steps.

\begin{center}
\renewcommand{\arraystretch}{1.1}
\begin{tabular}{|c|c|c|c|c|}
\hline
  Year & Reference & Symmetry & Global Asymmetry & Local Asymmetry \\
\hline\hline
  1992 & \cite{dacorogna-1} & $r=2$, $q\leq 2p$ & $r=2$, $q\gg 1$ &  \\
\hline
  1992 & \cite{dacorogna-1} & $r=q$ &  &  \\
\hline
  1997 & \cite{egorov} &  & $r=2$, $q>4p-1$ &  $r=2$, $q>4p-1$\\
\hline
  1998 & \cite{nazarov-0} &  & $q>3p$ &  \\
\hline
  1999--2000 & \cite{BK-1,BK-2} & $r=2$, $q\leq 2p+1$ &  &  \\
\hline
  2002 & \cite{nazarov-1} & $r=2$, $q\leq 3p$ &  &  \\
\hline
  2003 & \cite{dacorogna-2} & $q\leq rp+r-1$  & $q>(2r-1)p$ &  \\
\hline
  2011 & \cite{nazarov-2} & $q\leq(2r-1)p$ &  &  \\
\hline
  2017 & \cite{rovellini:tesi} &  &  & $r=2$, $q> 3p$ \\
\hline
\end{tabular}

\end{center}

From the technical point of view, the hardest step was proving symmetry of solutions to $\lambda_{W}(p,q,r)$ in the range $rp+r-1<q\leq(2r-1)p$. This is the content of~\cite{nazarov-1} in the case $r=2$ and of~\cite{nazarov-2} for general $r$. However, the proofs provided in these papers are ``quite technical'', in the sense that they require numerical computations carried out up to 18~significant digits in order to verify inequalities between functions with more than 10~levels of parentheses.

The contribution of this paper is twofold.
\begin{itemize}

  \item We provide a proof of the symmetry of minimizers to $\lambda_{W}(p,q,r)$ in the full range $q\leq(2r-1)p$ without computer assisted steps. Indeed, the key inequality of Proposition~\ref{prop:symmetry} is established through a suitable variable change, all whose details can be checked by a human.

  \item We prove nonexistence of odd \emph{local} minima to $\lambda_{W}(p,q,r)$ when $q>(2r-1)p$. This result was already known for \emph{global} minima in the same range, but for local minima it was limited to the case $r=2$ and $q>4p-1$ (but the argument extends in a straightforward way to $q>r^{2}p-(r-1)^{2}$ in the general case).

\end{itemize}

These two results settle completely the issue of symmetry/asymmetry of local and global minimizers to $\lambda_{W}(p,q,r)$.

This paper is organized as follows. In section~\ref{sec:literature} we review the main steps in previous literature that are needed in our approach. In section~\ref{sec:symmetry} we prove symmetry of minimizers for $q\leq(2r-1)p$. In section~\ref{sec:asymmetry} we prove asymmetry of local minimizers for $q>(2r-1)p$. 

\setcounter{equation}{0}
\section{Survey of previous literature}\label{sec:literature}

For the convenience of the reader, in this section we recall briefly the main approaches to Theorem~\ref{thmbibl:main} developed in the last decades. We focus in particular on the ideas that are needed in the sequel. 

\paragraph{\textmd{\emph{Poincaré inequality}}}

The symmetry of minimizers can be proved either via radial rearrangement (see~\cite{talenti}) or by inspecting the Euler equation associated to the variational problem. The Euler equation has a first integral, which up to rescaling and affine variable changes can be written in the form
\begin{equation}
|u'(x)|^{p}+|u(x)|^{q}=1.
\label{EE:P}
\end{equation}

From the theory of ordinary differential equations we obtain that all nonzero solutions to (\ref{EE:P}) are periodic with the same period, which we denote by $4T$. In particular, there exists a unique function $u_{P}:[-T,T]\to\re$ such that
\begin{eqnarray}
 & u_{P}'(x)=-\operatorname{sgn}(x)\left(1-|u_{P}(x)|^{q}\right)^{1/p}
\qquad
\forall x\in[-T,T],
 &  
 \label{defn:up-eqn}   \\[0.5ex]
 &  u_{P}(-T)=u_{P}(T)=0,  &	
\label{defn:up-dbc}  \\[0.5ex]
 &  u_{P}(x)>0
\qquad
\forall x\in(-T,T).  & 
\label{defn:up-ineq}
\end{eqnarray}

This is the ``unique'' equality case in Poincaré inequality, in the following sense.

\begin{propbibl}
Let $u:[a,b]\to\re$ be a local (in the $C^{1}$ norm) minimizer to (\ref{quotient}) among all  nontrivial functions $u\in C^{1}([a,b])$ with $u(a)=u(b)=0$. 

Then the following statements hold true.
\begin{enumerate}
\renewcommand{\labelenumi}{(\arabic{enumi})}

\item  There exists three positive real numbers $\alpha$, $\beta$, $\gamma$ such that
\begin{equation}
\alpha|u'(x)|^{p}+\beta|u(x)|^{q}=\gamma.
\qquad
\forall x\in[a,b].
\nonumber
\end{equation}

\item  If $u_{P}(x)$ denotes the solution to (\ref{defn:up-eqn})--(\ref{defn:up-ineq}), then
$$u(x)=u\left(\frac{a+b}{2}\right)u_{P}\left(\frac{2T}{b-a}(x-a)-T\right)
\qquad
\forall x\in[a,b].$$
\end{enumerate}
\end{propbibl}

This proves in particular that the graph of all minimizers, both local and global, is symmetric with respect to the vertical line through the middle point of $[a,b]$.

\paragraph{\textmd{\emph{Wirtinger inequality}}}

In this case the Euler equation, up to rescaling and affine variable changes, has a first integral of the form
\begin{equation}
|u'(x)|^{p}+|u(x)|^{q}=1+\mu|u(x)|^{r-2}u(x),
\label{EE:W}
\end{equation}
where the Lagrange multiplier $\mu$ comes from the integral constraint (\ref{hp:W-r}). This equation coincides with (\ref{EE:P}) when $\mu=0$. In particular, there exists a unique function $u_{W}:[-T,T]\to\re$ such that
\begin{eqnarray}
  & u_{W}'(x)=\left(1-|u_{W}(x)|^{q}\right)^{1/p}
\qquad
\forall x\in[-T,T], &  
\label{defn:uw-eqn}   \\[0.5ex]
  &  u_{W}'(-T)=u_{W}'(T)=0,  &  
\label{defn:uw-nbc}   \\[0.5ex]
  &    u_{W}'(x)>0
  \qquad
  \forall x\in(-T,T).
\label{defn:uw-monot}
\end{eqnarray}

Actually, the function $u_{W}(x)$ can be obtained from $u_{P}(x)$ with a cut-an-paste procedure analogous to (\ref{cut-and-paste}).

This function is the ``unique'' local and global minimum point to $\lambda_{W}(p,q,r)$ when \emph{restricted to odd functions} (namely functions whose graph is symmetric with respect to the middle point of $[a,b]$), and hence it is also the unique \emph{odd} candidate to be a minimum point without the symmetry condition.

\begin{propbibl}

Let us consider the minimization problem for the quotient (\ref{quotient}) with the integral constraint (\ref{hp:W-r}).

Then the following statements hold true.
\begin{enumerate}
\renewcommand{\labelenumi}{(\arabic{enumi})}

\item  If $u:[a,b]\to\re$ is any local (in the $C^{1}$ norm) minimizer, then there exists three positive real numbers $\alpha$, $\beta$, $\gamma$, and a real number $\delta$, such that
\begin{equation}
\alpha|u'(x)|^{p}+\beta|u(x)|^{q}=\gamma+\delta|u(x)|^{r-2}u(x)
\qquad
\forall x\in[a,b].
\nonumber
\end{equation}

Moreover, $u(x)$ is odd if and only if $\delta=0$.

\item  Let $u:[a,b]\to\re$ be a local (in the $C^{1}$ norm) minimizer in the class of odd functions. If $u_{W}(x)$ denotes the solution to (\ref{defn:uw-eqn})--(\ref{defn:uw-monot}), then
$$u(x)=\frac{b-a}{2T}u'\left(\frac{a+b}{2}\right)u_{W}\left(\frac{2T}{b-a}(x-a)-T\right)
\qquad
\forall x\in[a,b].$$
\end{enumerate}

\end{propbibl}

\paragraph{\textmd{\emph{An auxiliary function}}}

Let $M_{q,r}$ denote the set of real numbers $\mu$ such that the equation
\begin{equation}
1+\mu|x|^{r-2}x-|x|^{q}=0
\label{eqn:mqr}
\end{equation}
has at least two real solutions, and let $x_{1}(\mu)<0<x_{2}(\mu)$ denote the two solutions closest to the origin. The set $M_{q,r}$ is always a connected open set with center in the origin, and  when $q>r-1$ it turns out that $M_{q,r}=\re$ and (\ref{eqn:mqr}) has always exactly two solutions. We note that in any case $x_{1}(\mu)\in(-1,0)$ and $x_{2}(\mu)>1$ for every positive element $\mu\in M_{q,r}$.

Following \cite{nazarov-0}, let us consider the function $J_{p,q,r}:M_{q,r}\to\re$ defined by
\begin{equation}
J_{p,q,r}(\mu):=\int_{x_{1}(\mu)}^{x_{2}(\mu)}\left(1+\mu|x|^{r-2}x-|x|^{q}\right)^{1/p_{*}}\,dx
\qquad
\forall\mu\in M_{q,r}.
\label{defn:J}
\end{equation}

We observe that, from the geometric point of view, $J_{p,q,r}(\mu)$ represents one half of the area of the oval subset of the Euclidean plane defined as
$$\left\{(x,y)\in\re^{2}:|y|^{p_{*}}+|x|^{q}\leq 1+\mu|x|^{r-2}x\right\}.$$

The following result clarifies the deep connection between the minimization of the area of this oval and the minimization problem for $\lambda_{W}(p,q,r)$.

\begin{propbibl}[Connection between $J_{p,q,r}(\mu)$ and Wirtinger inequality]\label{propbibl:WJ}

Let $p>1$, $q>1$, $r>1$ be three real numbers, let $\lambda_{W}(p,q,r)$ be the optimal constant in Wirtinger inequality, and let $J_{p,q,r}(\mu)$ be defined by (\ref{defn:J}).

Then the following statements hold true.
\begin{enumerate}
\renewcommand{\labelenumi}{(\arabic{enumi})}

\item  It turns out that
\begin{equation}
\lambda_{W}(p,q,r)=\theta^{\theta}p_{*}^{1/p_{*}}q^{1/q}\min\left\{J_{p,q,r}(\mu):\mu\in M_{q,r}\right\}.
\nonumber
\end{equation}

\item  There exists an odd global minimizer for $\lambda_{W}(p,q,r)$ if and only if $\mu=0$ is a minimum point for $J_{p,q,r}(\mu)$.

\item  If there exists a local minimizer for $\lambda_{W}(p,q,r)$ which is not odd, then $J_{p,q,r}(\mu)$ admits a stationary point $\mu\neq 0$.

\end{enumerate}

\end{propbibl}

Proposition~\ref{propbibl:WJ} above reduces an infinite dimensional variational problem to the minimization of an integral function of just one real variable. 

\begin{propbibl}[Qualitative behavior of $J_{p,q,r}(\mu)$]\label{propbibl:J}

Let $p>1$, $q>1$, $r>1$ be three real numbers.

Then the function $J_{p,q,r}(\mu)$ defined in (\ref{defn:J}) is an even function of class $C^{2}$ in $\re$, and the following statements hold true.
\begin{enumerate}
\renewcommand{\labelenumi}{(\arabic{enumi})}

\item  In the range $q>(2r-1)p$ it turns out that
\begin{equation}
J_{p,q,r}''(0)<0.
\label{thbibl:J''}
\end{equation}

\item  In the range $q\leq(2r-1)p$ it turns out that
\begin{equation}
J_{p,q,r}'(\mu)>0
\quad\quad
\forall\mu>0.
\label{thbibl:J'}
\end{equation}
\end{enumerate}
\end{propbibl}

Combining Proposition~\ref{propbibl:WJ} and Proposition~\ref{propbibl:J} it follows that
\begin{itemize}
  \item when $q\leq(2r-1)p$ the unique stationary point of $J_{p,q,r}(\mu)$ is $\mu=0$, and therefore $u_{W}(x)$ is the ``unique'' global and local minimizer to $\lambda_{W}(p,q,r)$.

  \item when $q>(2r-1)p$ the minimum of $J_{p,q,r}(\mu)$ is not achieved for $\mu=0$, and therefore the \emph{global} minimizers to $\lambda_{W}(p,q,r)$ are not odd functions.
  
\end{itemize}

\paragraph{\textmd{\emph{Previous symmetry results}}}

All proofs of symmetry results are based on (\ref{thbibl:J'}). In order to obtain this inequality, the first step consists in writing the derivative as
\begin{equation}
J_{p,q,r}'(\mu)=\frac{1}{p_{*}}\int_{x_{1}(\mu)}^{x_{2}(\mu)}\frac{|x|^{r-2}x}{\left(1+\mu|x|^{r-2}x-|x|^{q}\right)^{1/p}}\,dx.
\nonumber
\end{equation} 

Then with an affine variable change one can transform the integrals in $[x_{1}(\mu),0]$ and $[0,x_{2}(\mu)]$ into integrals in $[0,1]$. Setting
$$m:=-\frac{x_{1}(\mu)}{x_{2}(\mu)},
\qquad\qquad
R:=\frac{1-m^{q}}{1+m^{r-1}},$$
with some algebra one finds that $J_{p,q,r}'(\mu)$ is equal, up to positive quantities, to
\begin{equation}
\int_{0}^{1}\left(\frac{x^{r-1}}{\left(1-R+Rx^{r-1}-x^{q}\right)^{1/p}}-\frac{m^{r}x^{r-1}}{\left(1-R-Rm^{r-1}x^{r-1}-m^{q}x^{q}\right)^{1/p}}\right)\,dx,
\label{key-integral}
\end{equation}
so that the symmetry result is equivalent to proving that this integral is positive for every $m\in(0,1)$ ($m$ lies in this interval because $|x_{2}(\mu)|>|x_{1}(\mu)|$ for positive $\mu$).

\begin{itemize}
  \item In the first paper~\cite{dacorogna-1} it was proved that the integrand is positive for every $(x,m)\in(0,1)^{2}$, and hence a fortiori also the integral is positive, when $r=2$ and $q\leq 2p$.
  
  \item In \cite{BK-1,BK-2} the previous argument was refined, and it was shown that the integrand is positive when $r=2$ and $q\leq 2p+1$.
    
  \item In~\cite{dacorogna-2} the refined argument was extended to general $r$, proving that the integrand is positive when $q\leq rp+r-1$.
  
  \item  It can be shown that when $q> rp+r-1$ the integrand is negative in a neighborhood of $x=1$, and hence the previous approach cannot be extended. Nevertheless, this does not imply that the integral is negative.
  
  \item  In~\cite{nazarov-1,nazarov-2} a computer assisted proof is provided in order to show that the integral is positive in the range $q\leq (2r-1)p$.
\end{itemize}

In Proposition~\ref{prop:symmetry} below we provide an elementary proof that (\ref{key-integral}) is positive for every $m\in(0,1)$ in the full range $q\leq(2r-1)p$.

\paragraph{\textmd{\emph{Previous asymmetry results}}}

We have seen that Proposition~\ref{propbibl:WJ} and Proposition~\ref{propbibl:J} imply that $u_{W}(x)$ is not a \emph{global} minimum point when $q>(2r-1)p$. The key tool is inequality (\ref{thbibl:J''}), obtained in \cite{nazarov-0} in the case $r=2$, and in \cite{nazarov-2} for general $r$. 

On the other hand, this is not enough to exclude that $u_{W}(x)$ remains a \emph{local} minimum point also for larger values of $q$. In previous literature the nonexistence of odd local minima was known just in a more restrictive range. The approach was the following (see~\cite{egorov}).

Let $[a_{0},b_{0}]\subseteq\re$ be an interval with $a_{0}=-b_{0}$, and let $u_{0}:[a_{0},b_{0}]\to\re$ be any odd function. For every $\ep\in(0,1)$ let us consider the function
\begin{equation}
u_{\ep}(x):=\left\{
\begin{array}{l@{\qquad}l}
\dfrac{1}{(1+\ep)^{1/(r-1)}}u_{0}\left(\dfrac{x}{1+\ep}\right) & \mbox{if }x\in[0,(1+\ep)b_{0}], \\[3ex]
\dfrac{1}{(1-\ep)^{1/(r-1)}}u_{0}\left(\dfrac{x}{1-\ep}\right) & \mbox{if }x\in[(1-\ep)a_{0},0].
\end{array}
\right.
\label{eqn:vc-bibl}
\end{equation}

We observe that $u_{\ep}$ is defined in the interval $[(1-\ep)a_{0},(1+\ep)b_{0}]$, whose length is again $b_{0}-a_{0}$ because $a_{0}+b_{0}=0$. Up to a translation, we can also assume that this interval is exactly $[a_{0},b_{0}]$, and up to this translation it turns out that $\uep\to u_{0}$ in the $C^{1}$ norm as $\ep\to 0^{+}$. Moreover, $u_{\ep}$ is obtained from $u_{0}$ through a piecewise affine variable change and
$$\int_{(1-\ep)a_{0}}^{(1+\ep)b_{0}}|\uep(x)|^{r-2}\uep(x)\,dx=\int_{a_{0}}^{b_{0}}|u_{0}(x)|^{r-2}u_{0}(x)\,dx,$$
so that $\uep(x)$ is a competitor for $\lambda_{W}(p,q,r)$ whenever $u_{0}(x)$ is.

With some standard calculations one finds that 
\begin{equation}
\frac{\|u_{\ep}'\|_{p}}{\|u_{\ep}\|_{q}}=\frac{\|u_{0}'\|_{p}}{\|u_{0}\|_{q}}\left(1-\frac{q-r^{2}p+(r-1)^{2}}{2(r-1)^{2}}\ep^{2}+o(\ep^{2})\right)
\qquad
\mbox{as }\ep\to 0^{+},
\nonumber
\end{equation}
which proves that any odd function, and hence in particular $u_{W}(x)$, is not a local minimizer in the range $q>r^{2}p-(r-1)^{2}$. 

More recently, the third author~\cite{rovellini:tesi} proved that $u_{W}(x)$ is not a local minimizer when $r=2$ and $q>3p$. This result settles the matter completely in the case $r=2$. The idea in~\cite{rovellini:tesi} was to define $\uep(x)$ as the increasing solution to (\ref{EE:W}) with $\mu=\ep$ that vanishes at the origin, and suitably modified near one of the endpoints of the maximal interval where it is increasing in order to fulfill the integral constraint (\ref{hp:W-r}). Unfortunately this approach, when extended to general $r$, does not seem to fill the full gap between $(2r-1)p$ and $r^{2}p-(r-1)^{2}$.

In this paper we go back to the original idea of modifying $u_{W}(x)$ through a variable change. The novelty is that the variable change we devise in section~\ref{sec:asymmetry} is nonlinear. It is less general than (\ref{eqn:vc-bibl}) in the sense that it does not apply to any odd function, but just to $u_{W}(x)$, which however is the unique candidate to be a local minimizer. On the other hand, with this variable change we show that $u_{W}(x)$ is not a local minimizer as soon as $q>(2r-1)p$, and therefore whenever it is not a global minimizer. This provides a proof of the asymmetry result, both for global and local minima, independent of the computation of the second derivative of $J_{p,q,r}(\mu)$.

\setcounter{equation}{0}
\section{Symmetry of local and global minimizers}\label{sec:symmetry}

In this section we prove that the integral (\ref{key-integral}) is positive when $q\leq(2r-1)p$. As we have seen, this implies the symmetry of minimizers in the same range.

The basic tool is the following elementary, but nevertheless powerful, result.

\begin{lemma}\label{lemma:trivial}

Let $f:(0,1)\to(0,+\infty)$ and $g:(0,1)\to(0,+\infty)$ be two continuous functions. Let $u:[0,1]\to[0,1]$ be an increasing function of class $C^{1}$ such that $u(0)=0$, $u(1)=1$, and
\begin{equation}
u'(x)f(u(x))>g(x)
\quad\quad
\forall x\in(0,1).
\label{hp:sup-sol}
\end{equation}

Then it turn out that
\begin{equation}
\int_{0}^{1}f(x)\,dx>\int_{0}^{1}g(x)\,dx.
\nonumber
\end{equation}

\end{lemma}

\paragraph{\textmd{\emph{Proof}}}

Thanks to assumption (\ref{hp:sup-sol}), with the variable change $x=u(t)$ we deduce that
$$\int_{0}^{1}f(x)\,dx=\int_{0}^{1}f(u(t))u'(t)\,dt>\int_{0}^{1}g(t)\,dt,$$
which completes the proof.\qed\medskip

In the following two results we prove the monotonicity of two real functions.

\begin{lemma}\label{lemma:monotone-phi}
For every pair of real numbers $0\leq a<b$, the function 
\begin{equation}
\varphi_{a,b}(x):=\frac{1-x^{a}}{1-x^{b}}
\quad\quad
\forall x\in(0,1)
\nonumber
\end{equation}
is nonincreasing (and decreasing if $a>0$).
\end{lemma}

\paragraph{\textmd{\emph{Proof}}}

The function $(b-a)x^{b}-bx^{b-a}+a$ is nonincreasing in $(0,1)$  (decreasing if $a>0$), and vanishes in $x=1$. It follows that
$$(b-a)x^{b}-bx^{b-a}+a\geq 0
\qquad
\forall x\in(0,1),$$
and hence
$$\varphi_{a,b}'(x)=-\frac{x^{a-1}}{(1-x^{b})^{2}}\cdot\left((b-a)x^{b}-bx^{b-a}+a\strut\right)\leq 0
\qquad
\forall x\in(0,1),$$
with strict inequalities when $a>0$. This implies the required monotonicity.\qed

\begin{lemma}\label{lemma:monotone-psi}
For every pair of real numbers $0<a<b$, the function 
\begin{equation}
\psi_{a,b}(x):=\frac{(1-x^{a})(1+x^{b})}{1-x^{a+b}}
\quad\quad
\forall x\in(0,1)
\nonumber
\end{equation}
is decreasing.
\end{lemma}

\paragraph{\textmd{\emph{Proof}}}

Since $0<a<b$, for every $x\in(0,1)$ it turns out that $b\log x<a\log x< 0$. Therefore, since the function $z\to z^{-1}\sinh z$ is decreasing in $(-\infty,0)$, it follows that
$$\frac{\sinh(a\log x)}{a\log x}<\frac{\sinh(b\log x)}{b\log x},$$
which is equivalent to saying that
$$b(x^{-a}-x^{a})<a(x^{-b}-x^{b})
\quad\quad
\forall x\in(0,1).$$

It follows that
$$\psi_{a,b}'(x)=\frac{x^{a+b-1}}{(1-x^{a+b})^{2}}\cdot\left(b(x^{-a}-x^{a})-a(x^{-b}-x^{b})\right)<0
\qquad
\forall x\in(0,1),$$	
which implies the required monotonicity.\qed
\medskip

We are now ready to state and prove the key inequality.

\begin{prop}\label{prop:symmetry}

Let $p>1$, $q>1$, $r>1$, and $m\in(0,1)$ be real numbers, and let us set
\begin{equation}
R:=\frac{1-m^{q}}{1+m^{r-1}}.
\label{defn:R}
\end{equation}

Let us assume that $q\leq(2r-1)p$.

Then it turns out that
\begin{equation}
\int_{0}^{1}\frac{x^{r-1}}{\left(1-R+Rx^{r-1}-x^{q}\right)^{1/p}}\,dx>
\int_{0}^{1}\frac{m^{r}x^{r-1}}{\left(1-R-Rm^{r-1}x^{r-1}-m^{q}x^{q}\right)^{1/p}}\,dx.
\label{th:main}
\end{equation}

\end{prop}

\paragraph{\textmd{\emph{Proof}}}

For every $m\in(0,1)$, let us set
\begin{equation}
D(x):=\left[1-(1-m^{r-1})x^{r-1}\right]^{1/(r-1)}
\qquad
\forall x\in[0,1],
\label{defn:D}
\end{equation}
and
\begin{equation}
u(x):=\frac{mx}{D(x)}
\qquad\quad
\forall x\in[0,1].
\label{defn:u}
\end{equation}

Simple computations show that
\begin{equation}
0<m< D(x)<1
\qquad\quad
\forall(x,m)\in(0,1)^{2}.
\label{est:m-D}
\end{equation}

Moreover, it turns out that $u(0)=0$, $u(1)=1$, and
\begin{equation}
u'(x)=\frac{m}{[D(x)]^{r}}
\qquad\quad
\forall x\in(0,1),
\label{eqn:u'}
\end{equation}
so that $u$ is an increasing function. Let $f(x)$ and $g(x)$ denote the integrands in the left-hand side and in the right-hand side of (\ref{th:main}), respectively. If we show that (\ref{hp:sup-sol}) is satisfied, then (\ref{th:main}) follows from Lemma~\ref{lemma:trivial}. Keeping (\ref{defn:u}) and (\ref{eqn:u'}) into account, inequality (\ref{hp:sup-sol}) becomes (in the sequel we simply write $D$ instead of $D(x)$)
\begin{equation}
1-R-R(mx)^{r-1}-(mx)^{q}> D^{(2r-1)p}
\left(1-R+R\frac{(mx)^{r-1}}{D^{r-1}}-\frac{(mx)^{q}}{D^{q}}\right),
\nonumber
\end{equation}
or equivalently
\begin{equation}
1-R>Rm^{r-1}x^{r-1}\cdot\frac{1+D^{(2r-1)p-(r-1)}}{1-D^{(2r-1)p}}+m^{q}x^{q}\cdot\frac{1-D^{(2r-1)p-q}}{1-D^{(2r-1)p}}.
\label{est:th2}
\end{equation}

Let us estimate the two summands in the right-hand side. For the first one, we deduce from (\ref{defn:D}) that
$$x^{r-1}=\frac{1-D^{r-1}}{1-m^{r-1}},$$
and then we observe that the assumptions of Lemma~\ref{lemma:monotone-psi} are satisfied with $a:=r-1$ and $b:=(2r-1)p-(r-1)$. Therefore, from (\ref{est:m-D}) it follows that
\begin{eqnarray}
m^{r-1}x^{r-1}\cdot\frac{1+D^{(2r-1)p-(r-1)}}{1-D^{(2r-1)p}}  &  =  &  \frac{m^{r-1}}{1-m^{r-1}}\cdot\frac{(1-D^{r-1})(1+D^{(2r-1)p-(r-1)})}{1-D^{(2r-1)p}}
\nonumber \\[1ex]
  & < & \frac{m^{r-1}}{1-m^{r-1}}\cdot\frac{(1-m^{r-1})(1+m^{(2r-1)p-(r-1)})}{1-m^{(2r-1)p}}.
  \nonumber
\end{eqnarray}

As for the second summand in (\ref{est:th2}), we observe that $x<1$ and the assumptions of Lemma~\ref{lemma:monotone-phi} are satisfied with $a:=(2r-1)p-q$ and $b:=(2r-1)p$ (this is the point where it is essential that $q\leq (2r-1)p$). Therefore, from (\ref{est:m-D}) it follows that
\begin{equation}
m^{q}x^{q}\cdot\frac{1-D^{(2r-1)p-q}}{1-D^{(2r-1)p}}\leq m^{q}\cdot\frac{1-m^{(2r-1)p-q}}{1-m^{(2r-1)p}}.
\nonumber
\end{equation}

From the last two estimates it follows that (\ref{est:th2}) is proved if we can show that
$$1-R\geq Rm^{r-1}\cdot\frac{1+m^{(2r-1)p-(r-1)}}{1-m^{(2r-1)p}}+m^{q}\cdot\frac{1-m^{(2r-1)p-q}}{1-m^{(2r-1)p}}.$$

Plugging (\ref{defn:R}) into this inequality, we discover that it is actually an equality.\qed

\setcounter{equation}{0}
\section{Asymmetry of local and global minimizers}\label{sec:asymmetry}

In this section we prove that odd functions are not local minimizers to $\lambda_{W}(p,q,r)$ when $q>(2r-1)p$. To this end, we can limit ourselves to showing that the function $u_{W}(x)$ defined by (\ref{defn:uw-eqn})--(\ref{defn:uw-monot}) is not a local minimizer. Indeed, we have seen that this is the unique local minimizer in the class of odd functions.

To begin with, we show some relations between integrals of powers of $u_{W}(x)$.

\begin{lemma}\label{lemma:integrals}

Let $u_{W}:[-T,T]\to\re$ be the solution to (\ref{defn:uw-eqn})--(\ref{defn:uw-monot}).

Then for every real number $s\geq 0$ it turns out that
\begin{eqnarray}
\int_{-T}^{T}|u_{W}(x)|^{2s}\,dx  &  =  &  \frac{(2s+1)p_{*}+q}{(2s+1)p_{*}}\int_{-T}^{T}|u_{W}(x)|^{q+2s}\,dx  
\nonumber\\[1ex]
  &  =  &  \frac{(2s+1)p_{*}+q}{q}\int_{-T}^{T}|u_{W}'(x)|^{p}\cdot|u_{W}(x)|^{2s}\,dx.
\label{th:int-pqs}
\end{eqnarray}
\end{lemma}

\paragraph{\textmd{\emph{Proof}}}

For the sake of shortness we simply write $u(x)$ instead of $u_{W}(x)$. From (\ref{defn:uw-eqn}) it follows that
$$|u'(x)|^{p}+|u(x)|^{q}=1
\qquad
\forall x\in[-T,T].$$

Multiplying both sides by $|u(x)|^{2s}$, and integrating in $[-T,T]$, we deduce that
\begin{equation}
\int_{-T}^{T}|u'(x)|^{p}\cdot|u(x)|^{2s}\,dx+\int_{-T}^{T}|u(x)|^{q+2s}\,dx=\int_{-T}^{T}|u(x)|^{2s}\,dx.
\label{eqn:int-1}
\end{equation}

On the other hand, from (\ref{defn:uw-eqn}) it follows also that
$$\int_{-T}^{T}|u|^{q+2s}\,dx=\int_{-T}^{T}|u|^{2s}u\cdot\frac{|u|^{q-2}u}{(1-|u|^{q})^{1/p}}u'\,dx=-\frac{p_{*}}{q}\int_{-T}^{T}|u|^{2s}u\left[\left(1-|u|^{q}\right)^{1/p_{*}}\right]'\,dx.$$

Now from (\ref{defn:uw-eqn}) and (\ref{defn:uw-nbc}) we deduce that $|u(-T)|=|u(T)|=1$, and hence when we integrate by parts we find that
\begin{eqnarray}
\int_{-T}^{T}|u(x)|^{q+2s}\,dx  &  =  &  \frac{(2s+1)p_{*}}{q}\int_{-T}^{T}|u(x)|^{2s}u'(x)\left(1-|u(x)|^{q}\right)^{1/p_{*}}\,dx
\nonumber \\[1ex]
  &  =  &  \frac{(2s+1)p_{*}}{q}\int_{-T}^{T}|u(x)|^{2s}\cdot|u'(x)|^{p}\,dx.
\label{eqn:int-2}
\end{eqnarray}

Combining (\ref{eqn:int-1}) and (\ref{eqn:int-2}) we obtain (\ref{th:int-pqs}).\qed

\subsection*{Proof of asymmetry}

\paragraph{\textmd{\emph{A general variable change}}}

Let $[a_{0},b_{0}]\subseteq\re$ be a symmetric interval with $a_{0}=-b_{0}$, and let $u_{0}:[a_{0},b_{0}]\to\re$ be an odd function of class $C^{1}$. Let $\varphi:[a_{0},b_{0}]\setminus\{0\}\to\re$ be a bounded function of class $C^{1}$. Let us assume that $\varphi$ is odd, and let $M$ denote the supremum of $|\varphi(x)|$ in $[a_{0},b_{0}]\setminus\{0\}$. For every $\ep\in(0,1/M)$, let us consider the function $\gep:[a_{0},b_{0}]\to\re$ defined as
\begin{equation}
\gep(y):=1+\ep\varphi(y)
\qquad
\forall y\in[a_{0},b_{0}]\setminus\{0\}.
\label{defn:gep}
\end{equation}

Let us consider the solution $\yep(x)$ to the problem
\begin{equation}
\yep'=\gep(\yep),
\qquad
\yep(0)=0.
\label{defn:yep}
\end{equation}

Since $\gep(y)$ is bounded from below by a positive constant, there exists a (non symmetric) interval $[\aep,\bep]$ such that $\yep(\aep)=a_{0}$ and $\yep(\bep)=b_{0}$, and in addition
$$\yep\in C^{0}([\aep,\bep])\cap C^{1}\left([\aep,\bep]\setminus\{0\}\right).$$

Since $u_{0}(0)=0$, there exists a unique continuous function $\uep:[\aep,\bep]\to\re$ such that
\begin{equation}
\uep(x):=\left[\yep'(x)\right]^{1/(r-1)}u_{0}(\yep(x))
\qquad
\forall x\in[\aep,\bep]\setminus\{0\}.
\label{defn:uep}
\end{equation}

With the variable change $y=\yep(x)$ we obtain that
\begin{eqnarray*}
\int_{\aep}^{\bep}|\uep(x)|^{r-2}\uep(x)\,dx  &  =  &  \int_{\aep}^{\bep}\yep'(x)\cdot|u_{0}(\yep(x))|^{r-2}u_{0}(\yep(x))\,dx  \\
  &  =  &  \int_{a_{0}}^{b_{0}}|u_{0}(y)|^{r-2}u_{0}(y)\,dy.
\end{eqnarray*}	

This means that, whenever $u_{0}(x)$ satisfies the integral constraint (\ref{hp:W-r}), the function $\uep(x)$ satisfies the same condition for every admissible value of $\ep$.

We can also transform $\uep(x)$ in a new function defined in $[a_{0},b_{0}]$ through a further affine variable change, and in this sense $\uep\to u_{0}$ in the $C^{1}$ norm as $\ep\to 0^{+}$.

\paragraph{\textmd{\emph{Expansion of the length of the interval}}}

We claim that the length of the interval satisfies
\begin{equation}
(\bep-\aep)^{\theta}=(b_{0}-a_{0})^{\theta}\left(1+\frac{I_{0}}{J_{0}}\theta\ep^{2}+o(\ep^{2})\right)
\qquad
\mbox{as }\ep\to 0^{+},
\label{taylor:theta}
\end{equation}
where
\begin{equation}
I_{0}:=\int_{a_{0}}^{b_{0}}\varphi^{2}(x)\,dx,
\qquad\qquad
J_{0}:=\int_{a_{0}}^{b_{0}}1\,dx=b_{0}-a_{0}.
\nonumber
\end{equation}

Indeed, from (\ref{defn:yep}) we obtain that
$$\bep-\aep=\int_{\aep}^{\bep}1\,dx=\int_{\aep}^{\bep}\frac{\yep'(x)}{\gep(\yep(x))}\,dx=
\int_{a_{0}}^{b_{0}}\frac{1}{\gep(y)}\,dy.$$

Moreover, from (\ref{defn:gep}) we deduce that
$$\frac{1}{\gep(y)}=\frac{1}{1+\ep\varphi(y)}=1-\ep\varphi(y)+\ep^{2}\varphi^{2}(y)+o(\ep^{2}).$$

Integrating in $[a_{0},b_{0}]$, and recalling that $\varphi(y)$ is an odd function, we conclude that
$$\bep-\aep=\int_{a_{0}}^{b_{0}}\frac{1}{\gep(y)}\,dy=(b_{0}-a_{0})+\ep^{2}\int_{a_{0}}^{b_{0}}\varphi^{2}(x)\,dx+o(\ep^{2}).$$

Raising both sides to the power $\theta$, we obtain exactly (\ref{taylor:theta}).

\paragraph{\textmd{\emph{Expansion of the norm of $\uep$}}}

We claim that the norm of $\uep$ in $L^{q}((\aep,\bep))$ satisfies
\begin{equation}
\|\uep\|_{q}=\|u_{0}\|_{q}\left(1+\frac{I_{1}}{J_{1}}\gamma_{1}\ep^{2}+o(\ep^{2})\right)
\qquad
\mbox{as }\ep\to 0^{+},
\label{taylor:q}
\end{equation}
where
$$\gamma_{1}:=\frac{q}{2(r-1)^{2}}-\frac{3}{2(r-1)}+\frac{1}{q},$$
$$I_{1}:=\int_{a_{0}}^{b_{0}}|u_{0}(x)|^{q}\cdot\varphi^{2}(x)\,dx,
\qquad\qquad
J_{1}:=\int_{a_{0}}^{b_{0}}|u_{0}(x)|^{q}\,dx.$$

Indeed, from (\ref{defn:yep}) we obtain that
\begin{eqnarray*}
\int_{\aep}^{\bep}|\uep(x)|^{q}\,dx  &  =  &  \int_{\aep}^{\bep}|u_{0}(\yep(y))|^{q}\cdot\left[\gep(\yep(x))\right]^{\frac{q}{r-1}-1}\cdot\yep'(x)\,dx   \\
  &  =  &  \int_{a_{0}}^{b_{0}}|u_{0}(y)|^{q}\cdot\left[\gep(y)\right]^{\frac{q}{r-1}-1}\,dy.
\end{eqnarray*}

Moreover, from (\ref{defn:gep}) we deduce that
$$\left[\gep(y)\right]^{\frac{q}{r-1}-1}=\left[1+\ep\varphi(y)\right]^{\frac{q}{r-1}-1}=1+\ep\left(\frac{q}{r-1}-1\right)\varphi(y)+q\gamma_{1}\ep^{2}\varphi^{2}(y)+o(\ep^{2}).$$

When we plug this expansion into the integral, the first order term cancels due to the symmetries, and we conclude that
$$\int_{\aep}^{\bep}|\uep(x)|^{q}\,dx=\int_{a_{0}}^{b_{0}}|u_{0}(y)|^{q}\,dy+q\gamma_{1}\ep^{2}\int_{a_{0}}^{b_{0}}|u_{0}(y)|^{q}\cdot\varphi^{2}(y)\,dy+o(\ep^{2}).$$

Raising both sides to the power $1/q$, we obtain exactly (\ref{taylor:q}).

\paragraph{\textmd{\emph{Expansion of the norm of $\uep'$}}}

Let us assume that there exists a continuous odd function $\psi:[a_{0},b_{0}]\to\re$ such that
\begin{equation}
\varphi'(x)u_{0}(x)=\psi(x)u_{0}'(x)
\qquad
\forall x\in[a_{0},b_{0}]\setminus\{0\}.
\label{defn:psi}
\end{equation}

Then we claim that the norm of $\uep'$ in $L^{p}((\aep,\bep))$ satisfies
\begin{equation}
\|\uep'\|_{p}=\|u_{0}'\|_{p}\left(1+\frac{I_{2,1}\gamma_{2,1}+I_{2,2}\gamma_{2,2}+I_{2,3}\gamma_{2,3}}{J_{2}}\ep^{2}+o(\ep^{2})\right)
\qquad
\mbox{as }\ep\to 0^{+},
\label{taylor:p}
\end{equation}
where
$$\gamma_{2,1}:=\frac{r^{2}p}{2(r-1)^{2}}-\frac{3}{2}\frac{r}{r-1}+\frac{1}{p},
\qquad
\gamma_{2,2}:=\frac{(p-1)}{2(r-1)^{2}},
\qquad
\gamma_{2,3}:=\frac{rp}{(r-1)^{2}}-\frac{2}{r-1},$$
and
$$J_{2}:=\int_{a_{0}}^{b_{0}}|u_{0}'(x)|^{p}\,dx,
\qquad\qquad
I_{2,1}:=\int_{a_{0}}^{b_{0}}|u_{0}'(x)|^{p}\cdot\varphi^{2}(x)\,dx,$$
$$I_{2,2}:=\int_{a_{0}}^{b_{0}}|u_{0}'(x)|^{p}\cdot\psi^{2}(x)\,dx,
\qquad\quad
I_{2,3}:=\int_{a_{0}}^{b_{0}}|u_{0}'(x)|^{p}\cdot\varphi(x)\cdot\psi(x)\,dx.$$

Indeed, from (\ref{defn:gep}) and (\ref{defn:yep}) we obtain that 
$$\yep''(x)=\gep'(\yep(y))\cdot\yep'(x)=\ep\varphi'(\yep(x))\cdot\yep'(x)
\qquad
\forall x\in[\aep,\bep]\setminus\{0\},$$
and hence, keeping (\ref{defn:psi}) into account, the time-derivative of (\ref{defn:uep}) turns out to be
$$\uep'(x)=\left[\yep'(x)\right]^{1/(r-1)}\cdot u_{0}'(\yep(x))\cdot\left(\frac{\ep\psi(\yep(x))}{r-1}+\gep(\yep(x))\right)
\qquad
\forall x\in[\aep,\bep]\setminus\{0\},$$
and therefore (for the sake of shortness, we do not write the explicit dependence on $x$ in the integrals of the first line)
\begin{eqnarray*}
\int_{\aep}^{\bep}|\uep'|^{p}\,dx  &  =  &  \int_{\aep}^{\bep}|u_{0}'(\yep)|^{p}\cdot\left|\frac{\ep\psi(\yep)}{r-1}+\gep(\yep)\right|^{p}\cdot\left[\gep(\yep)\right]^{\frac{p}{r-1}-1}\cdot\yep'\,dx   \\[1ex]
  &  =  &  \int_{a_{0}}^{b_{0}}|u_{0}'(y)|^{p}\cdot\left|\frac{\ep\psi(y)}{r-1}+\gep(y)\right|^{p}\cdot\left[\gep(y)\right]^{(p-r+1)/(r-1)}\,dy.
\end{eqnarray*}  

Recalling (\ref{defn:gep}), we find that the last integrand is equal to
\begin{eqnarray*}
\lefteqn{ |u_{0}'(y)|^{p}\cdot\left\{1+p\left(\varphi(y)+\frac{\psi(y)}{r-1}\right)\ep+\frac{p(p-1)}{2}\left(\varphi(y)+\frac{\psi(y)}{r-1}\right)^{2}\ep^{2}\right\}\cdot}
\\[1ex]
  &  &  \qquad\mbox{}\cdot\left\{1+\frac{p-r+1}{r-1}\varphi(y)\ep+\frac{(p-r+1)(p-2r+2)}{2(r-1)^{2}}\varphi^{2}(y)\ep^{2}\right\}+o(\ep^{2}).
\end{eqnarray*}

Now with some patience we multiply and we integrate in $[a_{0},b_{0}]$, and we observe that the terms of order~1 cancel due to the symmetries. Then we raise the result to the power $1/p$ and we get exactly (\ref{taylor:p}).

\paragraph{\textmd{\emph{Non-optimal asymmetry result -- Back to literature}}}

Just to present a somewhat trivial example, let us consider for a while the case where $\varphi(x)$ is the piecewise constant function equal to $-1$ in $(a_{0},0)$ and equal to $1$ in $(0,b_{0})$. Due to the lack of regularity of $\varphi$, the variable change $\yep(x)$ is just piecewise affine, and actually the function $\uep(x)$ defined in (\ref{defn:uep}) coincides with the one defined in (\ref{eqn:vc-bibl}). In this case (\ref{defn:psi}) holds true with $\psi(x)\equiv 0$, and the expansions (\ref{taylor:theta}), (\ref{taylor:q}), and (\ref{taylor:p}) hold true with
$$I_{0}=J_{0},
\qquad
I_{1}=J_{1},
\qquad
I_{2,1}=J_{2},
\qquad
I_{2,2}=I_{2,3}=0,$$
and therefore
$$(\bep-\aep)^{\theta}\frac{\|\uep'\|_{p}}{\|\uep\|_{q}}=(b_{0}-a_{0})^{\theta}\frac{\|u_{0}'\|_{p}}{\|u_{0}\|_{q}}\left(1+\Gamma_{p,q,r}\ep^{2}+o(\ep^{2})\right),$$
where
$$\Gamma_{p,q,r}  = \theta-\gamma_{1}+\gamma_{2,1}=-\frac{q-r^{2}p+(r-1)^{2}}{2(r-1)^{2}}.$$

This implies that every odd function $u_{0}(x)$ is not a local minimum point in the range $q>r^{2}p-(r-1)^{2}$, as already observed.

\paragraph{\textmd{\emph{Optimal asymmetry result}}}

Let us finally focus on the case where $u_{0}(x)=u_{W}(x)$, namely on the unique candidate to be a local minimum point in general. Let us choose $\varphi(x):=|u_{W}(x)|^{r-2}u_{W}(x)$, and let us observe that this function satisfies (\ref{defn:psi}) with $\psi(x):=(r-1)\varphi(x)$. As a consequence, the coefficients in (\ref{taylor:theta}) through (\ref{taylor:p}) can be computed by applying Lemma~\ref{lemma:integrals} with $s=0$ and $s=r-1$. We find that expansion (\ref{taylor:theta}) holds true with
$$I_{0}=\int_{a_{0}}^{b_{0}}|u_{W}(x)|^{2r-2}\,dx,
\qquad\qquad
J_{0}=b_{0}-a_{0},$$
expansion (\ref{taylor:q}) holds true with
$$I_{1}=\int_{a_{0}}^{b_{0}}|u_{W}(x)|^{q+2r-2}\,dx=\frac{(2r-1)p_{*}}{(2r-1)p_{*}+q}I_{0},
\qquad
J_{1}=\int_{a_{0}}^{b_{0}}|u_{W}(x)|^{q}\,dx=\frac{p_{*}}{p_{*}+q}J_{0},$$
and expansion (\ref{taylor:p}) holds true with
$$I_{2,1}=\frac{I_{2,2}}{(r-1)^{2}}=\frac{I_{2,3}}{r-1}=\int_{a_{0}}^{b_{0}}|u_{W}'(x)|^{q}\cdot|u_{W}(x)|^{2r-2}\,dx=\frac{q}{(2r-1)p_{*}+q}I_{0},$$
$$J_{2}=\int_{a_{0}}^{b_{0}}|u_{W}'(x)|^{p}\,dx=\frac{q}{p_{*}+q}J_{0}.$$

Putting things together, with some calculations we end up with
$$(\bep-\aep)^{\theta}\frac{\|\uep'\|_{p}}{\|\uep\|_{q}}=
(b_{0}-a_{0})^{\theta}\frac{\|u_{W}'\|_{p}}{\|u_{W}\|_{q}}\left(1+\frac{I_{0}}{J_{0}}\Gamma_{p,q,r}\ep^{2}+o(\ep^{2})\right),$$
where
\begin{eqnarray*}
\Gamma_{p,q,r}  &  =  &  \theta+\frac{p_{*}+q}{(2r-1)p_{*}+q}\left\{-(2r-1)\gamma_{1}+\gamma_{2,1}+(r-1)^{2}\gamma_{2,2}+(r-1)\gamma_{2,3}\right\}   \\[1ex]
  & = &  -\frac{2r-1}{2(r-1)^{2}}\cdot\frac{p_{*}+q}{(2r-1)p_{*}+q}\cdot\left(q-(2r-1)p\strut\right).
\end{eqnarray*}

This implies that $u_{W}(x)$ is not a local minimum point in the range $q>(2r-1)p$.\qed

%\bibliographystyle{MaxNew}
%\bibliography{Wirtinger}

\label{NumeroPagine}

\end{document}